\documentstyle[12pt]{article}

\begin{titlepage}

\title{Theta and zeta functions for
odd-dimensional locally symmetric
spaces of rank one}

\author{Ulrich
 Bunke\thanks{Humboldt-Universit\"at zu Berlin,
 Institut
 f\"ur Reine Mathematik (SFB288), Ziegelstr. 13a, Berlin 10099.
E-mail:ubunke@mathematik.hu-berlin.de
} and
Martin
Olbrich\thanks{
Humboldt-Universit\"at zu Berlin, Institut f\"ur Reine
Mathematik (SFB288), Ziegelstr.
13a, Berlin 10099.
E-mail:olbrich@mathematik.hu-berlin.de  }
}

\end{titlepage}

\begin{document}

\newcommand{\oh }{{\bf h}}
\newcommand{\om}{{\bf m}}

\newcommand{\kaaa}{{\bf k}}

\newcommand{\paaa}{{\bf p}}

\newcommand{\taaa}{{\bf t}}

\newcommand{\haaa}{{\bf h}}

\newcommand{\R}{{\bf R}}

\newcommand{\Z}{{\bf Z}}

\newcommand{\C}{{\bf C}}

\newcommand{\HR}{H{\bf R}}

\newcommand{\HC}{H{\bf C}}

\newcommand{\HH}{H{\bf H}}

\newcommand{\PR}{P{\bf R}}

\newcommand{\PC}{P{\bf C}}

\newcommand{\PH}{P{\bf H}}
\newcommand{\laaa}{{\bf l}}
\newcommand{\Saaa}{{\bf S}}

\newcommand{\G}{{\bf G}}

\newcommand{\A}{{\bf A}}

\newcommand{\Naaa}{{\bf N}}

\newcommand{\g}{{\bf g}}

\newcommand{\aaa}{{\bf a}}

\newcommand{\naaa}{{\bf n}}

\newcommand{\M}{{\bf M}}

\newcommand{\K}{{\bf K}}

\newcommand{\Oaaa}{{\cal O}}

\newcommand{\Haaa}{{\bf H}}

\newcommand{\db}{{\bar{\partial}}}

\newcommand{\Paaa}{{\bf P}}

\maketitle

\newtheorem{prop}{Proposition}[section]

\newtheorem{lem}[prop]{Lemma}

\newtheorem{ddd}[prop]{Definition}

\newtheorem{theorem}[prop]{Theorem}

\newtheorem{kor}[prop]{Corollary}

\newtheorem{ass}[prop]{Assumption}

\newtheorem{con}[prop]{Conjecture}

\tableofcontents
\section{Introduction}

This paper is a continuation of our work on theta and zeta functions
\cite{bunkeolbrich93},
\cite{bunkeolbrich94} and \cite{bunkeolbrich941}.
In the previous papers
we considered the case
of even dimensional
rank one symmetric spaces
of non-compact type. The present is concerned
with the odd-dimensional case, i.e. with odd-dimensional real hyperbolic
manifolds.

It is the natural appearence of eta invariants in connection with super
theta and zeta functions what makes this case particularly interesting.
That the eta invariant is connected with special value of a zeta function
was first observed for the signature operator by Millson \cite{millson78}
and later by Moscovici/Stanton \cite{moscovicistanton89} in the general case
(even for higher rank spaces). The super theta function was introduced by
Juhl \cite{juhl93}. He derived its properties using the super zeta function
while
the latter is studied through dynamical Lefschetz formulas
and by specializing the results of \cite{moscovicistanton89}.

In odd dimensions the analytic torsion (introduced by Ray/Singer
\cite{raysinger71})
is another interesting spectral
invariant. It is associated to a finite dimensional unitary representation
of the fundamental group yielding a flat vector bundle over the locally
symmetric space that is used to twist the Laplacians on differential forms.
The analytic torsion is then given by a certain combination of
zeta-regularized determinants of these twisted Laplacians. It  was first
observed by Fried \cite{fried862}
that the value of the twisted Ruelle zeta function at zero is the inverse of
the
analytic torsion. The connection of special values of zeta functions with
analytic torsion
and (higher) analytic torsion was generalized to higher rank cases by
Moscovici/Stanton \cite{moscovicistanton91}
and Deitmar \cite{deitmar94}.

Thus, the odd-dimensional case is was in some sense exhausted.
Moreover, in this case the identity contribution to the Selberg trace
formula is much simpler than in the even dimensional case since the the
Plancherel densities are polynomials.
In other terms, the distributional local trace of the corresponding
wave operators is a finite sum of derivatives of delta distributions
concentrated at zero (see Corollary \ref{fr56}),
i.e. a weak sort of Huygens' principle holds.  This a priori
excludes the difficulties with normalization constants one was confronted
with in the even dimensional case and which were a major motivation for us to
write
\cite{bunkeolbrich93},\cite{bunkeolbrich94}. There we found an easy solution
using the connection of the analysis on the symmetric spaces and its compact
dual space.
Though not as essential as in the even dimensional case   we apply this duality
the present paper  in order
to give an easy derivation of the identity contribution to the Selberg trace
formula
avoiding harmonic analysis on the non-compact symmetric space.

Our main motivation to write this paper was to show that one can recover all
previous results on theta and zeta functions by a "unique continuation" of our
methods of \cite{bunkeolbrich94}. We believe that our results, even if not
explicitly stated  in the literature, are known to the specialists in the
field.
But we think that our approach is rather short and much less involved
then the previous ones.

As in \cite{bunkeolbrich94} our plan is as  follows.
Let $\G=SO(n,1)$ or $\G=Spin(n,1)$ and $\K=SO(n)$  or $\K=Spin(n)$
be a maximal compact subgroup.
The zeta and theta functions
are associated to $\M$-types, $\M$
 being defined as the centralizer of an
 Iwasawa $\aaa$.
We first obtain $\K$-types restricting
 to a given $\M$-type.
We use this $\K$-types in order to define associated vector bundles over
the locally symmetric space and certain differential operators.
We then prove trace formulas for functions of these operators.
Once having the trace formula we study the theta, the Selberg zeta and
the Ruelle zeta function in a way similar to \cite{bunkeolbrich94}.
Because of the application to analytic torsion we include twists,
i.e. finite dimensional unitary representations  of the  fundamental group
of the locally symmetric space. Like in the present paper it is  easy to
include twists in
\cite{bunkeolbrich94} as well.
\section{The trace formula}
\subsection{The restriction map $r:R(\K)\rightarrow R(\M)$}

We consider the pairs of groups $(\K,\M)$ with $\K=Spin(n)$ and
$\M=Spin(n-1)$ or $\K=SO(n)$ and $\M=SO(n-1)$ for an odd integer
$n\geq 3$. Let $R(\K)$ and $R(\M)$ be the corresponding representation
rings over $\Z$. In this subsection we shall investigate the restriction map
\[ r: R(\K)\longrightarrow R(\M)\]
induced by the inclusion $\M\hookrightarrow\K$.
Let $\g=\kaaa\oplus\aaa\oplus\naaa$ be an Iwasawa decomposition of
 $\g=so(1,n)$ such that
$\M\subset \K$ can be identified with the centralizer of $\aaa$.
Let $W:=W(\g,\aaa)$ be the Weyl group of the restricted root system and
$w\in W$ its non-trivial element. $w$ acts on $\hat \M$ by
\[ w\sigma(m):= \sigma(m_w^{-1}m m_w)\,\qquad m\in \M, \sigma\in \hat\M,\]
where $m_w$ is a representative of $w$ belonging to the normalizer
of $\aaa$ in $\K$. This action extends to $R(\M)$. Then $R(\M)$
splits (over $\Z[\frac{1}{2}] $) into the
 $\pm 1$-eigenspaces of $w$
\[  R(\M)\otimes \Z[\frac{1}{2}]  =R(\M)^+\otimes \Z[\frac{1}{2}]\oplus
R(\M)^-\otimes\Z[\frac{1}{2}]   \ .\]

We consider the Cartan
 subalgebra $\taaa$ of $\om= {\bf so}(n-1)$, which can also be
considered via the inclusion $\om\hookrightarrow \kaaa$ as a Cartan
subalgebra of $\kaaa= {\bf so}(n)$, given by
\[ \taaa := \left\{ T_\mu:= \left( \begin{array}{ccc}
\begin{array}{cc}
0&-\mu_1\\ \mu_1&0
\end{array}&&\\
&\ddots&\\
&&\begin{array}{cc}
0&-\mu_k\\ \mu_k&0
\end{array}
\end{array}\right)\:|\:\mu_i\in\R, k=\frac{n-1}{2}\right\}\ .\]
We denote by $\nu=(\nu_1,\dots,\nu_k)$ the element in $\imath\taaa^*$,
which sends $T_\mu$ to $\imath( \nu_1\mu_1+\dots+\nu_k\mu_k)$.
 We
 choose systems of positive roots
$$\Phi^+(\kaaa^c,\taaa):= \{e_i\pm e_j (i<j),e_i\}\ \mbox{and}\
\Phi^+(\om^c,\taaa):= \{e_i\pm e_j (i<j)\}\  ,$$
where $e_i:=(0,\dots ,0,1,0,\dots,0)$ with $1$ at the $i$'th entry.
Their Weyl groups are given by
\begin{eqnarray*}
W_k&=&\{\mbox{permutations of the
 coordinates with possible sign changes}\}\ ,\\
W_m&=&\{\mbox{permutations of the
coordinates with an even number of sign changes}\}\ ,
\end{eqnarray*}
and we have
\begin{eqnarray*} \rho_k:=\frac{1}{2}\sum_{\alpha\in\Phi^+(
\kaaa^c,\taaa)}\alpha &=&(k-\frac{1}{2},k-\frac{3}{2},\dots,
\frac{3}{2},\frac{1}{2})\ ,\\
\rho_m:=\frac{1}{2}\sum_{\alpha\in\Phi^+(\om^c,\taaa)}\alpha
&=&(k-1,k-2,\dots,1,0)\ .
\end{eqnarray*}
Now by the theory of the highest weights
 we can write
\[ \hat \K = \{\gamma_\nu\:|\:\nu_1\geq\dots
\geq\nu_k\geq 0,
\nu_i\in \Z ,i=1,\dots ,k\ \mbox{or, if $\K=
Spin(n)$,}\ \nu_i\in \frac{1}{2}\Z\}\]
and
\[ \hat \M = \{\sigma_\nu\:|\:\nu_1\geq\dots
\geq\nu_{k-1}\geq |\nu_k|, \nu_i\in \Z ,i=1,\dots ,k\ \mbox{or,
 if $\M=Spin(n-1)$,}\ \nu_i\in \frac{1}{2}\Z\}\ .\]
Then the action of $w$ on $\hat\M$ looks as follows :
$ w\sigma_{(\nu_1,\dots,\nu_k)}=\sigma_{(\nu_1,\dots,-\nu_k)}$.
In particular, for the half-spin representations $s^\pm=\sigma_{(1/2,\dots,\pm
1/2)}$ of $Spin(n-1)$ we have $ws^\pm=s^\mp$. We denote by $s$ the
spin representation of $Spin(n)$. The first part of the following
proposition can already be found in Miatello/Vargas \cite{miatellovargas83}.
\begin{prop}\label{martin}
The map $r$ has the following properties :
\begin{enumerate}
\item $r$ is a bijection between $R(\K)$ and $R(\M)^+$.
\item If $\sigma\in R(\M)^-$, then there exists a unique element
$\gamma\in R(Spin(n))$ with $\sigma=(s^+-s^-)r(\gamma)$ such
that $s\cdot\gamma\in R(\K)$.
\item More explicitly, if $\nu_k>0$, then
\begin{equation}\label{qui} \sigma_{(\nu_1,\dots,\nu_k)}-w
\sigma_{(\nu_1,\dots,\nu_k)} =(s^+-s^-)r(\gamma_{(\nu_1-
\frac{1}{2},\dots,\nu_k-\frac{1}{2})})\ .
\end{equation}
In addition, the $\K$-representation $s\otimes
\gamma_{\nu-(\frac{1}{2},\dots,\frac{1}{2})}$ splits into
two representations $\gamma^+(\nu)$ and $\gamma^-(\nu)$
 such that
\begin{equation}\label{qqui} \sigma_{(\nu_1,\dots,\nu_k)}+w
\sigma_{(\nu_1,\dots,\nu_k)} =r(\gamma^+(\nu)-\gamma^-(\nu))\ .
\end{equation}
\end{enumerate}
\end{prop}
{\it Proof:} We employ the Weyl character formula. Let us denote
 by $\chi_\sigma$ ($\chi_\gamma$) the character of an irreducible
 representation $\sigma$ of $Spin(n-1)$ ($\gamma$ of $Spin(n)$).
We introduce the Weyl denominators
\begin{eqnarray*}
\Delta_{\K}(expT)&:=& \prod_{\alpha\in \Phi^+(\kaaa^c,\taaa)}
(e^{\alpha(T)/2} - e^{-\alpha(T)/2})\ ,\quad T\in\taaa,\\
\Delta_{\M}(expT)&:=& \prod_{\alpha\in \Phi^+(\om^c,\taaa)}
(e^{\alpha(T)/2} - e^{-\alpha(T)/2})\ .
\end{eqnarray*}
Note that $\chi_{s^+}(expT)-\chi_{s^-}(expT)=\Delta_{\K}(expT)/
\Delta_{\M}(expT)$. For a highest weight $\nu$ of $Spin(n)$ we compute
\begin{eqnarray}\label{mo}
&&((\chi_{s^+}-\chi_{s^-})\chi_{\gamma_\nu})(expT)\nonumber\\
&=&(\chi_{s^+}-\chi_{s^-})(expT)\frac{\sum_{s\in W_k} dets
e^{s(\nu+\rho_k)(T)}}{\Delta_\K(expT)}\nonumber\\
&=&\frac{\sum_{s\in W_k} dets
e^{s(\nu+\rho_k)(T)}}{\Delta_\M(expT)}\nonumber\\
&=&\frac{\sum_{s\in W_k} dets
e^{s(\nu+(1/2,\dots,1/2)+\rho_m)(T)}}{\Delta_\M(expT)}\nonumber\\
&=&\frac{\sum_{s\in W_m} dets
(e^{s(\nu+(1/2,\dots,1/2)+\rho_m)(T)}
-e^{s(\nu+(1/2,\dots,1/2,-2\nu_k-1/2)+\rho_m)(T)})}{\Delta_\M(expT)}
\nonumber\\&=& \chi_{\sigma_{\nu+(1/2,\dots,1/2)}}(expT)
-\chi_{w\sigma_{\nu+(1/2,\dots,1/2)}}(expT)\ .
\end{eqnarray}
Thus we have proved (\ref{qui}), and the second part of the proposition
follows. Further we conclude that $r$ is injective, since in view of (\ref{mo})
 the map $R(K)\ni \gamma\rightarrow (s^+-s^-)r(\gamma)\in R(Spin(n-1))$
 has no kernel. In fact
$$\sigma_{\nu+(1/2,\dots,1/2)}\not=w\sigma_{\nu+(1/2,\dots,1/2)}.$$
Since the action of $w$ interchanges
the $\M$-isotypic components of an
 element in $\hat \K$, we see that the image
 of $r$ is contained in $R(\M)^+$.
In order to prove the first assertion of the
 proposition it remains to construct
 a preimage in R(\K) for every $\sigma\in
R(\M)^+$. Since $(s^+-s^-)\sigma\in R(Spin(n-1))^-$,
we find an element $\gamma\in R(Spin(n))$,
which actually belongs to R(\K), such
 that
\[ (s^+-s^-)\sigma=(s^+-s^-)r(\gamma)\ .\]
Now $R(Spin(n-1))$ is a polynomial ring,
 hence $\sigma=r(\gamma)$.

We are left with the proof of (\ref{qqui}).
 Let us consider the character of $s\otimes\gamma_{\nu}$
\begin{eqnarray}
\chi_s\chi_{\gamma_\nu}(expT)&=&
\sum_{v\in\{\pm1\}^k}\sum_{s\in W_k} \frac{detse^{(s(\nu+\rho_k)+v(\frac{1}{2},
\dots,\frac{1}{2}))(T)}}{\Delta_\K(expT)}\nonumber\\
&=& \sum_{v\in\{\pm1\}^k}\sum_{s\in W_k}
\frac{detse^{s(\nu+v(\frac{1}{2},\dots,\frac{1}{2})+\rho_k)(T)}} {
\Delta_\K(expT)}\nonumber\\
&=& \sum_{\{v\in\{\pm1\}^k\:|\:\nu+v(\frac{1}{2},\dots,\frac{1}{2})
+\rho_k \mbox{\scriptsize regular in } \taaa^*\}}\sum_{s\in W_k}
\frac{detse^{s(\nu+v(\frac{1}{2},\dots,\frac{1}{2})+\rho_k)(T)}}{
\Delta_\K(expT)}\ . \label{split}
\end{eqnarray}
In the last step we have used that for singular $\xi$
\[ \sum_{s\in W_k}
detse^{s\xi(T)}\equiv 0\ .\]
Observe that $\nu+v(\frac{1}{2},\dots,\frac{1}{2})$ is a
highest weight of an irreducible representation of $\K$ iff
$\nu+v(\frac{1}{2},\dots,\frac{1}{2})+\rho_k$ is regular.
Therefore we can define
\[ \gamma^\pm(\nu):= \sum_{\{v\in\{\pm1\}^k\:|\:\nu+v(
\frac{1}{2},\dots,\frac{1}{2})+\rho_k \mbox{\scriptsize
regular in } \taaa^*, detv=\pm 1\}} \gamma_{\nu+v(
\frac{1}{2},\dots,\frac{1}{2})}\ ,\]
and (\ref{split}) says that $s\otimes
\gamma_\nu=\gamma^+(\nu)\oplus\gamma^-(\nu)$.

For $\xi=(\xi_1,\dots,\xi_k)$ we set $\xi^-:=(\xi_1,\dots,\xi_{k-1},-\xi_k)$.
 Let $p$ be the projection $p:W_k\rightarrow \{\pm 1\}^k$ given by
$p(s)((\frac{1}{2},\dots,\frac{1}{2})):=s((\frac{1}{2},\dots,\frac{1}{2}))$.
Then we have
\begin{eqnarray*}
&&(\Delta_\K(
\chi_{\sigma_{\nu+(1/2,\dots,1/2)}} +\chi_{w\sigma_{\nu+
(1/2,\dots,1/2)}}))(expT)\\
&=&(\Delta_\M(\chi_{s^+}-\chi_{s^-}) (\chi_{\sigma_{\nu+
(1/2,\dots,1/2)}} +\chi_{w\sigma_{\nu+(1/2,\dots,1/2)}}))(expT)\\
&=&\left(\sum_{v\in\{\pm1\}^k}detve^{v(\frac{1}{2},\dots,
\frac{1}{2})(T)}\right) \left(\sum_{s\in W_m}
dets(e^{s(\nu+\rho_k)(T)}+e^{s((\nu+
\rho_k)^-)(T)})\right)\\
&=&\left(\sum_{v\in\{\pm1\}^k}
detve^{v(\frac{1}{2},\dots,\frac{1}{2})(T)}\right)\left(\sum_{s\in W_k} dets\
det p(s)e^{s(\nu+\rho_k)(T)}\right)\\
&=&\sum_{v\in\{\pm1\}^k}\sum_{s\in W_k} det v\  det s
e^{s(\nu+v(\frac{1}{2},\dots,\frac{1}{2})+\rho_k)(T)}\\
&=&(\Delta_\K(\chi_{\gamma^+(\nu)} -\chi_{\gamma^-(\nu)}))(expT)\ .
\end{eqnarray*}
Now (\ref{qqui}) follows. $\Box$

\subsection{The contribution of the identity}
Let $\sigma\in\hat{\M}$.
We have to distinguish the two cases (a): $\sigma=w\sigma$ and (b):
$\sigma\not=w\sigma$, where $w$ is the non-trivial element of the Weyl group of
$(\g,\aaa)$.
In the following we apply Proposition \ref{martin} several times.
Let $\gamma\in R(\K)$ be the unique lift of $\sigma$, i.e
$r(\gamma)=\sigma$ in the case (a) and $r(\gamma)=\sigma+w\sigma$ in the case
(b).
In the case (b) we also consider
 the super lift $\gamma^s=s\otimes \gamma^\prime\in R(\K)$.
$s$ is the spinor representation
 of $Spin(n)$ and $\gamma^\prime\in Spin(n)^\wedge$ is irreducible  such
that $\sigma-w\sigma=\pm(s^+-s^-)
\otimes r(\gamma^\prime)$ holds in $R(Spin(n-1))$, where
$s^\pm$ are the half spinor representations of $Spin(n-1)$
and  $r:R(\K)\rightarrow R(\M)$ (or $r:R(Spin(n))\rightarrow R(Spin(n))$)
is the restriction homomorphism.
If $\K,\M\not=Spin(n),Spin(n-1)$ we view $R(\K),R(\M)$ as subrings of
$R(Spin(n)),R(Spin(n-1))$.
Let $\g=\kaaa\oplus\paaa$ be the Cartan decomposition of $\g$.
Then there is a Clifford multiplication $\paaa\otimes s\rightarrow s$.
In fact, there are two Clifford multiplications differing by a sign.
Let $\aaa\subset\paaa$ be the one-dimensional subspace and $H\in\aaa$ be
the unit vector pointing into the positive direction of $\aaa$ given by the
Iwasawa decomposition.
Then we fix the sign of the Clifford multiplication (and hence of the
Dirac operator later on) by requireing that $H$ acts as multiplication by
$\pm \imath$ on $s^\pm\subset s$. Since $\M\subset \K$ is defined as the
centralizer
of $\aaa$ and thus $H$ is $\M$-invariant the subspaces $s^\pm\subset s$ are
well defined and $H$ respects the decomposition $s=s^+\oplus s^-$.
Note that by the third part of Proposition \ref{martin}
we can choose the representative of the lift $\gamma$ of
$\sigma +w\sigma$ such that $|\gamma|=|\gamma^s|$, i.e. both representations
coincide up to the $Z_2$-grading.
Here we abuse the notation denoting the element of $R(\K)$ and its
representative,
i.e. a formal sum of irreducible representations with integer coefficients by
the same
symbol $\gamma$. $|\gamma|$ is only defined for the representative.

Let $V^{(s)}=V(\gamma^{(s)})$ be the associated bundle over $X:=\G/\K$.
It is a $\G$-homogeneous bundle. Let $-\Omega$ be the Casimir operator
of $\G$ acting on sections of $V^{(s)}$ and $\Omega_\M$ be the Casimir operator
of $\M$. We normalize the invariant scalar product on $\g$ such that it induces
a metric of sectional curvature $-1$ on $X$. The invariant scalar product
of $\om$ is obtained by restriction. This fixes the normalization
of the Casimir operators. Let $\rho:=
\frac{n-1}{2} $ and $c(\sigma):=\rho^2 -\sigma(\Omega_\M)$.
We define the operator $(A^{(s)})^2:=\Omega-c(\sigma)$.
It gives rise to an unbounded selfadjoint operator on $L^2(X,V^{(s)})$.
Let $A^{(s)}$ be its square root with non-negative real and imaginary part.
In the super case $V^s$ is a Dirac bundle carrying
the Dirac operator $D$ such that $A^s=|D|$.

We define the distributions $K(t):= tr\:cos(tA)(x,x)$
and $J(t):=tr\: D\:cos(tA^s)(x,x)$, where $x\in X$ is arbitrary.
Applying Hadamard's analysis as in \cite{bunkeolbrich94}
we obtain the structure of these distributions at $t=0$.

We first consider $K(t)$.
Exactly as in \cite{bunkeolbrich94}, Corollary 3.2
  we derive  applying \cite{hoermander85}, Thm. 17.5.5.
\begin{lem}\label{zumverg}
The distribution $t\rightarrow K(t)$ has an asymptotic
 expansion at $t\to 0$ of the form
$$K(t)=\sum_{k=-(n-1)/2}^0 c_{2k}  \delta^{(2k)}(t)+
 \sum_{k=0}^\infty c_{2k+1} |t|^{2k+1}.$$
\end{lem}
In fact, we will show that $c_{k}=0$ for $k>0$ employing the analysis
on the dual symmetric space $S^n=X_d=\G^d/\K$, i.e. the unit sphere.
Let $V_d\rightarrow X_d$ be the homogeneous bundle associated to the lift
$\gamma$ carrying the operator $A_d$.
Let $\g=\kaaa\oplus\aaa\oplus\naaa$ be  the Iwasawa decomposition of $\g$,
$\alpha$ be the root of $(\aaa,\naaa)$ and $H_\alpha\in\aaa$ be the root
vector satisfying $\alpha(H_\alpha)=1$.
Define $\epsilon(\sigma)\in\{0,1/2\}$ by the condition :
$e^{2\pi\imath\epsilon(\sigma)}=\sigma(2\pi\imath H_\alpha)\in\{\pm 1\}$.
We consider the the lattice $L:=\epsilon(\sigma)+\Z$ and the
even Weyl polynomial $P(\lambda):=P(\lambda,\sigma)$ (\cite{bunkeolbrich94}).
In order to define $P(\lambda,\sigma)$ let
$\om\subset\kaaa$ be the centralizer of $\aaa$
 and $\oh=\taaa\oplus\aaa \subset\om\oplus\aaa$
be a Cartan subalgebra of $\g$.
Then $\taaa$ is a Cartan subalgebra of $\om$.
Fix a positive root system $\Phi^+(\g,\oh)$
compatible with the choice of $\alpha$, let
$\delta:=\frac{1}{2}\sum_{\beta\in\Phi^+(\g,\oh)}\beta$ and
$\rho_m:=\delta-\frac{n-1}{2} \alpha$.
Then $$P(\lambda,\sigma):=\prod_{\beta\in\Phi^+(\g,\oh)}
\frac{(\lambda+\mu_\sigma+\rho_m,\beta)}{(\delta,\beta)},$$
where $\mu_\sigma\in \imath\taaa^\ast$ is the highest
 weight of $\sigma$.
Note that $\epsilon(w\sigma)=\epsilon(\sigma)$
and $P(\lambda,w\sigma)=P(\lambda,\sigma)$.

We consider weighted dimensions of eigenspaces since
$\gamma$ is a virtual representation.
Let $m_d(\lambda,\gamma,\sigma):=TrE_{A_d}
(\{\lambda\})$ be the multiplicity of the eigenvalue $\lambda\in\R$.
Moreover, for $0\le \lambda\in L$ let
$m_d(\lambda):=m_d(\lambda,\sigma):=P(\lambda,\sigma)$ in the case (a)
and $m_d(\lambda):=m_d(\lambda,\sigma):=2 P(\lambda,\sigma)$ in the case (b).
{}From \cite{bunkeolbrich94}, Proposition 9.3
we obtain
\begin{lem} \label{ppp1}
Essentially, the eigenvalues of $A_d$ form the ladder
$\{0\le \lambda\in L\}$ and have the multiplicities
$m_d(\lambda,\gamma,\sigma)=m_d(\lambda)$.
 There are at most finitely many
exceptional eigenvalues or gaps of finite multiplicity.
\end{lem}
Now we compute the distribution $K_d(t):=Tr\: cos(tA_d)$.
Note the equation
$$\sum_{k=0}^\infty cos(tk)=\pi\delta(t)+1/2.$$
Thus, if $\epsilon(\sigma)=0$, then
$$\tilde{K}_d(t):=d\:\sum_{k=0}^\infty P(k) cos(tk)=d
 \pi P(\imath \frac{d}{dt}) \delta(t)+d P(0)/2$$
(note that $P(\lambda)$ is even), where $d=1$
in the case (a) and $d=2$ in the case (b).
Also
\begin{eqnarray*}
\sum_{k=0}^\infty cos((k+1/2)t)&=&cos(t/2)
\sum_{k=0}^\infty cos(kt)-sin(t/2)\sum_{k=0}^\infty sin(kt)\\
&=&\pi\delta(t)+\frac{1}{2}cos(t/2)-\frac{1}{2}cos(t/2)\\
&=&\pi\delta(t)
\end{eqnarray*}
Hence in the case $\epsilon(\sigma)=1/2$
$$\tilde{K}_d(t):=d\:\sum_{k=0}^\infty P(k+1/2)
 cos((k+1/2)t)=d \pi P(\imath \frac{d}{dt}) \delta(t).$$

Applying Hadamard's analysis \cite{hoermander85}, Thm. 17.5.5
to the wave operator $cos(tA_d)$
we get an asymptotic expansion
$$K_d(t)=\omega_{n+1}\sum_{k=-(n-1)/2}^0
c_{d,2k}  \delta^{(2k)}(t)+ \sum_{k=0}^\infty c_{d,2k+1} |t|^{2k+1},$$
where $\omega_{n+1}:=vol(S^n)$.
Observe that by Proposition \ref{ppp1}
$$K_d(t)-\tilde{K}_d(t)=\sum_{r=0}^{\mbox{finite}} c_r cos(t\lambda_r),$$
where $c_r,\lambda_r\in \C$.
Comparing this result with Lemma \ref{zumverg}  we obtain:
\begin{prop}
We have $$K_d(t)=\pi P(\imath \frac{d}{dt}) \delta(t)$$
and $m_d(\lambda,\sigma,\gamma)=m_d(\lambda)$ for all $0\not=\lambda\in L$
and in the case $\epsilon(\sigma)=0$ we have $m_d(0,\sigma,\gamma)=d P(0)/2$.
\end{prop}
Thus the uniqe lift $\gamma\in R(\K)$ of $\sigma\in\hat{\M}$
is admissible in the sense of \cite{bunkeolbrich94}, Definition  2.2.

By \cite{bunkeolbrich94}, Lemma 3.6, we have
$c_{d,k}=\imath c_k$.
Thus $c_k=0$ for $k>0$ and we obtain
\begin{kor}\label{fr56}
The distribution $K(t)$ is given by
$$K(t)=d\pi\omega_{n+1}^{-1} P(\frac{d}{dt}) \delta(t).$$
\end{kor}

\begin{lem}
The distribution $J(t)$ vanishes identically.
\end{lem}
{\it Proof:}
In analogy to the derivation of \ref{fr56} we must show the vanishing
of the corresponding distribution $J_d(t)$ on the compact dual space.
But the non-zero spectrum of $D_d$ is symmetric, i.e.
$dim\:E_{D_d}(\{\lambda\})= dim\:E_{D_d}(\{-\lambda\})$,
$\forall\lambda\not=0$.
$\Box$\newline

\subsection{The hyperbolic contribution}
The hyperbolic contribution is defined for semisimple $g\in \G$
to be the distribution
\begin{eqnarray*}
C_c^\infty(\R)\ni\phi\rightarrow I_\phi(g)&:=&\int_{\Omega_{<g>}} tr
K_\phi(x,gx) dx\\
C_c^\infty(\R)\ni\phi\rightarrow I^s_\phi(g)&:=&\int_{\Omega_{<g>}} tr
K^s_\phi(x,gx) dx,
\end{eqnarray*}
where $\Omega_{<g>}$ is the fundamental domain in $X$ of the cyclic group
generated by $g$ and $K_\phi^{(s)}(x,y)$ is the kernel
of the smoothing, finite propagation operator
\begin{eqnarray*}K_\phi&:=&\int_\R\phi(t)cos(tA)  dt\\
                 K^s_\phi&:=&\int_\R\phi(t)D \:cos(tA^s)  dt .
\end{eqnarray*}
We identify the fibres $V^{(s)}_x, V^{(s)}_{gx}$ by $g$
using the homogeneous bundle structure of $V^{(s)}$.

Let $\G=\K\A\Naaa$ be an adapted Iwasawa decomposition
of $\G$ such that $g=ma\in\M\A^+$. Let $l(g)=|log(a)|$ and
$$C(g,\sigma):=-\frac{l(g)e^{\rho l(g)} tr\:\sigma(m)}{2det(1-Ad(ma)_\naaa)}.$$
The proof of the following Proposition is similar to the one of Theorem 3.11 in
\cite{bunkeolbrich94}.
\begin{prop}
The hyperbolic contribution $I_\phi(g)$ is given by
\begin{eqnarray*}
I_\phi(g)&=&C(g,\sigma) (\phi(l(g)) +\phi(-l(g)))\quad\mbox{case (a)}\\
I_\phi(g)&=&(C(g,\sigma)+C(g,w\sigma))(\phi(l(g))+\phi(-l(g)))\quad\mbox{case
(b)}.
\end{eqnarray*}
\end{prop}
\begin{theorem}\label{flat}
The super hyperbolic contribution $I^s_\phi(g)$ is given by
$$I^s_\phi(g)=\imath(C(g,\sigma)-C(g,w\sigma))
(\phi^\prime(l(g))-\phi^\prime(-l(g))),$$
where $\phi^\prime(t)$ denotes the derivative of $\phi$ with respect to $t$.
\end{theorem}
{\it Proof:}
We first identify the domain $\Omega_{<g>}$.
We write $\G=\A\Naaa\K$ and $X=\A\Naaa$.
Then it is easy to see that one can choose $\Omega_{<g>}:=[1,a]\times
 \Naaa\subset \A\Naaa$.
Moreover, we trivialize
the bundle
\begin{equation}\label{trtr}
V^s =\A\Naaa\times V_{\gamma^s}
\end{equation}
identifying $V_{\gamma^s}=V_o^s$, $o=[\K]\in X$.
The kernel $K^s_\phi(x,y)$ becomes an $End(V_{\gamma^s})$-valued function.

We first carry out the integration with respect to $\A$.
 \begin{eqnarray*}
I^s_\phi(g)&=&\int_0^a\int_\Naaa tr\left(
K^s_\phi(bn,ab\alpha_m(n))\gamma^s(m) \right) db dn\\
&=&|\log(a)|\int_\Naaa tr\left( K^s_\phi(n,a\alpha_m(n))\gamma^s(m) \right)
dn\\
&=&l(g)\int_\Naaa tr\left( K^s_\phi(\alpha_{am}(n^{-1})n,a)\gamma^s(m) \right)
dn.\\
\end{eqnarray*}
The map $h(ma):n\rightarrow \alpha_{am}(n^{-1})n$ is a diffeomorphism of
$\Naaa$ and
the Haar measure $dn$ transforms as ( compare Helgason \cite{helgason84} 1.5.4
)
$$d h(ma)(n)=-det(1-Ad(ma)_\naaa)dn.$$
  We   continue the evaluation of $I^s_\phi(g)$ obtaining
\begin{equation}\label{zhnj}I^s_\phi(g)=
-\frac{l(g)}{det(1-Ad(ma)_\naaa)}
\int_\Naaa tr\left( K^s_\phi(a^{-1}n,o)
\gamma^s(m)\right) dn.\end{equation}

In order to evaluate the $\Naaa$-integral we employ the Fourier transform
of Helgason-type
$$F:C_c^\infty(X,V^s)\rightarrow
C^\infty(\imath\aaa^\ast\times\K,V_{\gamma^s}),$$
which is defined by
$$<v,F(f)(\lambda,k)>:=\int_X <f(x), \phi_{-\lambda,v}(k^{-1}x)>_x dx.$$
Here $v\in V_{\gamma^s}^\ast$ and
$\phi_{\lambda,v}\in C^\infty(X,V^s)$ is defined by
$\phi_{\lambda,v}(an):=a^{\lambda+\rho} v$ with respect to the trivialization
(\ref{trtr}).
The Fourier transform obviously extends to
 distributions with compact support $C_c^{-\infty}(X,V^s)$
$$F(f)(\lambda,k):=<f, l_k\phi_{-\lambda,v}>,$$
where $l_k$ is the left action by $k\in\K$.
Let $-\Omega$ and $\Omega_\M$ be the Casimir operators of $\G$ and $\M$. For
$\lambda\in\imath\aaa^\ast$  we have
$$\Omega\phi_{\lambda,v}=\phi_{\lambda,(|\lambda|^2+|\rho|^2-\Omega_\M)v}.$$
This can be seen by using the decomposition
$\g=\om\oplus \aaa  \oplus( \naaa+\bar{\naaa})$.
For $f\in C_c^\infty(X,V^s)$ we have
\begin{eqnarray*}
F((A^s)^2 f)&=&(|\lambda|^2+\rho^2-\Omega_\M-c(\sigma))F(f)\\
F(Df)&=&(\sum_{i=1}^{n-1}c(X_i)\gamma^\prime(Y_i)+\lambda c(H))F(f).
\end{eqnarray*}
We explain the notation of the second line.
We choose an orthonormal basis  $(H,X_i,i=1,\dots,n-1)$ of $\paaa$ with
$H\in\A$,
where $\g=\kaaa\oplus\paaa$ is the Cartan decomposition,
$Y_i\in \kaaa$, $X_i+Y_i\in \naaa$.
Then $c(X_i)$ is the Clifford multiplication by $X_i$ acting on the $s$-factor
of $\gamma^s$.

As in the proof of \cite{bunkeolbrich94}, Thm. 3.11 we show that
\begin{equation}\label{trew}F(D\:cos(tA^s) f)
=(\sum_{i=1}^{n-1}c(X_i)\gamma^\prime(Y_i)+\lambda c(H))
cos(t \sqrt{|\lambda|^2+\rho^2-\Omega_\M-c(\sigma)})F(f).\end{equation}

Let $I$ label the $\M$-types $\sigma^\prime$ occurring in $V_{\gamma^s}$.
We decompose
$$V_{\gamma^s}= \oplus_{\sigma^\prime\in I}V_\gamma(\sigma^\prime)$$
and let $P_{\sigma^\prime}$ be the corresponding projections.

By the finite propagation speed of $D \:cos(tA^s)$ and since $\phi$ has
compact support we have
$$K^s_\phi:C_c^\infty(X,V^s)\rightarrow C_c^\infty(X,V^s).$$
Integrating (\ref{trew}) against $\phi(t)$ we obtain
$$
F(K_\phi^s f)= (\sum_{i=1}^{n-1}
c(X_i)\gamma^\prime(Y_i)+\lambda c(H))
     \sum_{\sigma^\prime \in I}
\frac{1}{2}(\hat{\phi}(d(\lambda,\sigma^\prime))
+\hat{\phi}(-d(\lambda,\sigma^\prime)))  P_{\sigma^\prime} F(f),$$
where $$d(\lambda,\sigma^\prime)
:=\sqrt{|\lambda|^2+\rho^2-\Omega_\M-c(\sigma)}$$ and
$$\hat{\phi}(s)=\int_\R \phi(t)
 e^{\imath ts} dt.$$
This equation extends to
distributions $f\in C_c^{-\infty}(X,V^s)$.
Inserting for $f$ the delta
distribution located in $o$ with values in $V_{\gamma^s}^\ast$,
we obtain
\begin{eqnarray*}&&F(K^s_\phi(.,o))(\lambda,1)\\
&=&(\sum_{i=1}^{n-1}c(X_i)\gamma^\prime(Y_i)+\lambda c(H))
    \sum_{\sigma^\prime \in I} \frac{1}{2}(\hat{\phi}(d(\lambda,
\sigma^\prime))+\hat{\phi}(
-d(\lambda,\sigma^\prime)))
 P_{\sigma^\prime}\in End(V_\gamma^s).\end{eqnarray*}
Now we can continue the evaluation of the hyperbolic contribution.
Note that \linebreak[4]$K^s_\phi(a^{-1}n,o) \gamma^s(m)$ is a smooth
$End(V_{\gamma^s})$-valued
function with compact support.
We compute
\begin{eqnarray*}
&&a^{-\rho}\int_N  tr  K^s_\phi(a^{-1}n,o) \gamma^s(m)\\
&=&\frac{1}{2\pi\imath }\int_{\imath \aaa^\ast}
\int_\A e^{\lambda(log(a)+log(b))} b^{\rho} \int_N  tr
 K^s_\phi(bn,o)\gamma^s(m) dn\:db\:d\lambda\\
&=&\frac{1}{2\pi\imath}
 \int_{\imath \aaa^\ast}
 e^{\lambda(log(a))}tr F(K^s_\phi(.,o))
(-\lambda,1)\gamma^s(m) d\lambda\\
&=&\frac{1}{2\pi\imath}
\int_{\imath \aaa^\ast} e^{\lambda(log(a))}    tr
(\sum_{i=1}^{n-1}c(X_i)\gamma^\prime(Y_i)-
\lambda c(H))
 \\&&\quad \sum_{\sigma^\prime \in I}
 \frac{1}{2}(\hat{\phi}(d(\lambda,\sigma^\prime))+\hat{\phi}
(-d(\lambda,\sigma^\prime))) P_{\sigma^\prime}
   \gamma^s(m)d\lambda
\end{eqnarray*}
Note that $r(\gamma^s)=(s^+\otimes
 r(\gamma^\prime))\oplus
(s^-\otimes r(\gamma^\prime))$
as an $\M$-module and $\sum_{i=1}^{n-1}
c(X_i)\gamma^\prime
(Y_i)$ intertwines
these two summands.
Thus this sum does not contribute to the trace.
$c(H)$ acts by $\pm \imath$ on $s^\pm$.
Using $(s^+-s^-)\otimes r(\gamma^\prime)=\sigma-w\sigma$
we obtain further
\begin{eqnarray*}
&&a^{-\rho}\int_N  tr  K^s_\phi(a^{-1}n,o) \gamma^s\\
&=&-\frac{1}{2\pi} \int_{\imath \aaa^\ast}  \lambda e^{\lambda(log(a))}
\frac{1}{2}(\hat{\phi}(-|\lambda|)+\hat{\phi}(|\lambda|))
tr(\sigma(m)-w\sigma(m))  d\lambda\\
&=&\frac{\imath}{2} (\phi^\prime(l(g))-\phi^\prime(-l(g)))
tr (\sigma(m)-w\sigma(m)),
\end{eqnarray*}
where $\phi^\prime(t)=\frac{d}{dt}\phi(t)$.
The contribution of the   $\M$-types different from $\sigma$ and
$w\sigma$  cancels out
and $\sigma,w\sigma$ contribute with multiplicity $1,-1$.
We obtain the expression for $I^s_\phi(g)$ claimed in the Theorem
by inserting the last line into (\ref{zhnj}). $\Box$\newline

\subsection{The distributional trace formula}

Let $M=\Gamma\backslash\G/\K$ be an odd-dimensional,  closed,
hyperbolic manifold.
Let $\chi:\Gamma\rightarrow End(E)$ be a finite-dimensional unitary
representation of $\Gamma$ and $E_M$ be the corresponding flat bundle
over $M$. Let $r:=dim(E)$.
Let $(A^{(s)}_M)^2,D_M$ be the operators
associated to $\sigma\in\hat{\M}$ and twisted with $E$
acting on sections of  $V_M\otimes E_M\rightarrow M$.
Define $A^{(s)}_M$ as the square root of the selfadjoint extension of
$(A^{(s)}_M)^2$ with non-negative real and imaginary part.
We will suppress the $\chi$ in our notation, where it is possible.
Let $\phi\in C_c^\infty(\R)$. Then
\begin{eqnarray*}
K_{M,\phi}&:=&\int_{-\infty}^\infty \phi(t) cos(tA_M) dt\\
K^s_{M,\phi}&:=&\int_{-\infty}^\infty \phi(t)D_M\: cos(tA^s_M) dt
\end{eqnarray*}
are of trace class. The linear maps
$$C_c^\infty(\R)\ni\phi\rightarrow Tr K^{(s)}_{M,\phi}$$
define   distributions on $\R$ formally written as
$$Tr cos(tA_M),\quad\quad\quad Tr D_M\:cos(tA^s_M).$$

Let $n_\Gamma(g)$ be the number of classes in $\Gamma_g/<g>$, where
$\Gamma_g$ is the centralizer of $g$ in
 $\Gamma$ and $<g>$ is the
 group generated by $g$.
By $C\Gamma$ we denote the set of conjugacy classes of $\Gamma$.

\begin{prop}
We have the following equations of distributions.
\begin{eqnarray*}
&&Tr cos(tA_M)\\&=&vol(M)\pi r \omega^{-1}_{n-1}P(\frac{d}{dt})\delta(t)\\
&+&\sum_{[g]\in C\Gamma,[g]\not=[1]}
\frac{C(g,\sigma)tr\:\chi(g)}{n_\Gamma(g)}(
\delta(t-l(g))+\delta(-t-l(g)))\quad \mbox{case (a)},      \\[4mm]
&&\\&&Tr cos(tA_M)\\&=&  vol(M)2\pi r
\omega_{n+1}^{-1} P(\frac{d}{dt})\delta(t)\\
&+&\sum_{[g]\in C\Gamma,[g]\not=[1]} \frac{(C(g,\sigma)+
C(g,w\sigma))
tr\:\chi(g)}{
n_\Gamma(g)}(\delta(t-l(g))+\delta(-t-l(g)))\quad \mbox{case (b)},      \\[4mm]
&&\\&&Tr D\:cos(tA_M^s)\\
&=&\sum_{[g]\in C\Gamma,[g]\not=[1]}-\imath \frac{(C(g,\sigma)-C(g,w\sigma))
tr\:\chi(g)}{
n_\Gamma(g)}(\delta^\prime(t-l(g))-\delta^\prime(t+l(g)))
\end{eqnarray*}
\end{prop}
The proof is analogous to the one of \cite{bunkeolbrich94}, Prop 3.12.
The twist $\chi$ is incorporated by noting that
$$K_{M,\phi}(x,x)=\sum_{g\in\Gamma} K_\phi(\tilde{x},g\tilde{x})\otimes
\chi(g),$$
where $\tilde{x}\in X$ is a preimage of $x\in M$.
In the super case the contribution of the identity vanishes.
\section{Theta functions}

Let $M=\Gamma\backslash \G/\K$  be a closed, odd-dimensional
hyperbolic manifold, $\sigma\in \hat{\M}$ and $\chi$ be a finite dimensional
unitary representation of $\Gamma$.
Let $\gamma\in R(\K)$ be the lift of $\sigma$ and $(A^{(s)}_M)^2$ be
the associated operator on $M$ twisted with $\chi$.
Define $A_M^{(s)}$ as the square root of $(A^{(s)}_M)^2$
with non-negative real and imaginary part. In the case (b) let $D_M$ be
the corresponding Dirac operator twisted with $\chi$.
\begin{ddd}
For $Re(t)>0$ we define the theta function associated to $M$,
$\sigma$ and $\chi$ by
$$\theta(t):=\theta_\chi(t,\sigma):=Tr e^{-tA_M}.$$
In the super case (b) we also define the super theta function
$$\theta^s(t):=\theta^s_\chi(t,\sigma):=Tr sign(D_M) e^{-tA_M^s}.$$
\end{ddd}
The theta functions are  holomorphic functions on the right half plan
$Re(t)>0$. The goal of the present section is to provide their
meromorphic
extensions and to discuss their singularity picture.
We will suppress the $\sigma$ and $\chi$ in our notation.
\begin{theorem}
The theta function
 $\theta(t)$ admits a meromorphic continuation
 to the whole
complex plane. It satisfies the functional equation
$$\theta(t)+\theta(-t)=0.$$
The singularities of $\theta(t)$ are
first order poles at
$t=\pm \imath l(g)$, $g\in C\Gamma$ with residue
$$
res_{t=\pm \imath
l(g)}\theta(t)=\left\{\begin{array}{cc}\frac{C(g,\sigma)tr\:\chi(g)}{\pi
n_\Gamma(g)}&\mbox{case (a)}\\

\frac{(C(g,\sigma)+C(g,w\sigma))
tr\:\chi(g)}{\pi n_\Gamma(g)}&\mbox{case (b)}\end{array}\right.$$
and a pole of order $n-1$ at $t=0$.
The super theta function has a meromorphic
continuation to the whole
 complex plane.
It is regular at $t=0$ and satisfies
$$\theta^s(t)-\theta^s(-t)=0,\quad \theta(0)=\eta(D),$$
where $\eta(D)$ is the eta invariant of the Dirac operator $D$.
The singularities of $\theta^s(t)$ are
first order poles at
$t=\pm \imath l(g)$, $g\in C\Gamma$ with residue
$$
res_{t=\pm \imath l(g)}\theta(t)=\pm\frac{(C(g,\sigma)-C(g,w\sigma))
tr\:\chi(g)}{\pi n_\Gamma(g)}.$$
\end{theorem}
The properties of
the super theta function were already obtained by
Juhl  \cite{juhl93},Satz 9.1.4.
{\it Proof:}
The proof  for the non-super theta functions goes exactly as
in \cite{bunkeolbrich94}, Thm. 4.6. The only modification is
that we define $\theta$ on the half plane $Re(t)<0$ using the new functional
equation, i.e. $\theta(-t):=-\theta(t),Re(t)>0$.
Note that the identity contribution of the trace formula $K(t)$
is localized at $t=0$. The singularity of $\theta$ at $t=0$ is again given
by the local pseudodifferential analysis as in \cite{duistermaatguillemin75}.

We now discuss the super case. We first derive a meromorphic continuation
of the derivative $(\theta^s)^\prime(t)$. It turns out that this
derivative has second order poles with vanishing residues and can thus be
integrated to give $\theta^s$.

Since $D_M=sign(D_M)A_M^s$ we have
$(\theta^s)^\prime(t)=-Tr D_Me^{-tA_M^s}$.
We try to continue $(\theta^s)^\prime(t)$ using
 the functional
 equation
\begin{equation}\label{meckern}(\theta^s)^\prime(t)+(\theta^s)^\prime
(-t)=0.\end{equation}
We define $(\theta^s)^\prime$ as a distribution on $\C$.
Let $\phi\in C_c^\infty(\C)$. Then by definition
$$<(\theta^s)^\prime,\phi>:=\lim_{\epsilon\to 0}
\int_{|Re(v)|\ge \epsilon} (\theta^s)^\prime(v)\phi(v) dv.$$
We compute the distributional derivative $\bar{\partial}(\theta^s)^\prime$.
Note that $(\theta^s)^\prime(\imath t+\epsilon)$ converges
(considered as a distribution with respect to $t$)
to a tempered distribution on $\R$ when $\epsilon \to 0$.
\begin{eqnarray*}
&&<\bar{\partial}(\theta^s)^\prime,\phi>\\&=&<(
\theta^s)^\prime,\bar{\partial}^\ast \phi>\\
&=&lim_{\epsilon\to 0}\int_\epsilon^\infty\int_{
-\infty}^\infty \frac{1}{2}[-\frac{\partial}{\partial
u}-\frac{\partial}{\partial(\imath t)}]
(\phi(\imath t+
u)+\phi(-\imath t-u))\:\:(\theta^s)^\prime(\imath t+u)du\:\:dt\\
&=&<\frac{1}{2}((\theta^s)^\prime(\imath t+0)+(
\theta^s)^\prime(-\imath t+0)),\phi(\imath t) >.
\end{eqnarray*}
The distributional trace formula states
\begin{eqnarray*}
&&(\theta^s)^\prime(\imath t+0)+(\theta^s)^\prime(-\imath t+
0)\\&=&-2 Tr\: D_M\:cos(t A^s_M)\\
&=&2\imath \sum_{[g]\in C\Gamma,[g]\not=[1]}
\frac{(C(g,\sigma)-C(g,w\sigma))tr\:\chi(g)}{n_\Gamma(g)}(
\delta^\prime(t-l(g))-\delta^\prime(t+l(g))).
\end{eqnarray*}
Thus,
$$ <\bar{\partial}(\theta^s)^\prime,\phi>=-\imath\sum_{[g]\in
C\Gamma,[g]\not=[1]} \frac{(C(g,\sigma)-C(g,w\sigma))tr\:\chi(g)}{n_\Gamma(g)}(
\frac{d}{dt}_{|t=l(g)}\phi(\imath t) -\frac{d}{dt}_{|t=-l(g)}\phi(\imath t)).$$
Write $z=u+\imath t$.
Since
$$\bar{\partial}\frac{-\imath}{z^2}=\frac{\partial}{\partial t}(\frac{
\partial}{\partial u}-\frac{\partial}{\partial(\imath t)})\frac{1}{z} =\pi
\frac{\partial}{\partial
t}\delta(z)$$ we obtain that
$(\theta^s)^\prime(t)$ is holomorphic near the imaginary axis except
at $t=\pm \imath l(g)$, $g\in C\Gamma$, where it has second order
poles with the singular part
$$\pm\frac{(C(g,\sigma)-C(g,w\sigma))tr\:\chi(g)}{\pi n_\Gamma(g)}
\frac{1}{(t\mp \imath l(g))^2}.$$
Thus we have provided a meromorphic continuation of $(\theta^s)^\prime$
such that it only has second order poles with vanishing residues.
Note that $(\theta^s)^\prime(t)\to 0$ exponentially as $t\to\infty$.
Since the same holds for $\theta^s$ itself, we have
$$\theta^s(t):=\int_\infty^t (\theta^s)^\prime(u) du.$$
Thus $\theta^s$ is meromorphic and has first order poles
with the residues claimed in the statement of the theorem.
Integrating (\ref{meckern}) we obtain
$$
\theta^s(t)-\theta^s(-t)=0.
$$
Moreover,
\begin{equation}\label{lastline}
\theta^s(0)=-\int_0^\infty (\theta^s)^\prime(u) du
\end{equation}
We now have to identify the expression (\ref{lastline}) with $\eta(D)$.
The $\eta$-invariant of $D$ is defined as the value at $s=0$
of the $\eta$-function
$$\eta(s)=Tr\frac{sign(D)}{|D|^s}.$$
Initially, the $\eta$-function is defined for $Re(s)$ large. But
it has a meromorphic
continuation to all of $\C$ and turns out to be regular at $s=0$
(Atiyah/Patodi/Singer \cite{atiyahpatodisinger75}). We have
$$\eta(s)=Tr\frac{sign(D)}{|D|^s}=\frac{1}{\Gamma(s+1)}
\int_0^\infty u^s Tr\: D\:e^{-u|D|}du.$$
Since $-Tr \:D\:e^{-u|D|}=(\theta^s)^\prime(u)$ is regular at $u=0$
we can let $s=0$ obtaining
$$-\int_{0}^\infty(\theta^s)^\prime(u)du=\eta(D).$$
$\Box$\newline
\section{Zeta functions}
\subsection{The logarithmic derivative}
Let $M=\Gamma\backslash \G/\K$  be a closed, odd-dimensional
hyperbolic manifold, $\sigma\in \hat{\M}$ and $\chi$ be a finite dimensional
unitary representation of $\Gamma$.
We say that $[g]\in C\Gamma$ is primitive, if $n_\Gamma(g)=1$.
\begin{ddd}
For $Re(s)>\rho=(n-1)/2$ we define  the
Selberg zeta function associated to $\sigma$ by the Euler product.
$$Z_{S,\chi}(s,\sigma)=\prod_{[g]\in C
\Gamma,[g]\not=1,{\rm primitive}}\prod_{k=0}^\infty
\:det\left(1-e^{(-\rho-s)l(g)}S^k(Ad(g)_\naaa^{-1})\otimes
\sigma(m)\otimes \chi(g)\right).$$
\end{ddd}
We will again suppress $\sigma$ and $\chi$ in our notation,
where it is possible.
In the case (b) we study $Z_S(s,\sigma)$ via
$S(s):=S(s,\sigma):=Z_S(s,\sigma)Z_S(s,w\sigma)$ and
$S^s(s):=S^s(s,\sigma):=Z_S(s,\sigma)/Z_S(s,w\sigma)$.
$S^{(s)}$ are closely related to the analysis of the operator $D_M$.

In the case (a) the logarithmic derivative
of $Z_S(s)$ is given by (see e.g. \cite{bunkeolbrich94}, Prop. 6.3)
$$D(p):=\sum_{[g]\in C\Gamma,[g]
\not=1}\frac{2C(g,\sigma)tr \:\chi(g) e^{-pl(g)} }{n_\Gamma(g)}.$$
In the case  (b) the logarithmic
derivatives of $S^{(s)}(s)$ are given by
\begin{eqnarray*}
D(p)&:=&\sum_{[g]\in C\Gamma,[g]\not=1}2\frac{
(C(g,\sigma)+
C(g,w\sigma))tr \:\chi(g)e^{-pl(g)} }{n_\Gamma(g)}\\
D^s(p)&:=&\sum_{[g]\in C\Gamma,[g]\not=1}2\frac{(C(
g,\sigma)-C(g,w\sigma))
tr\: \chi(g)e^{-pl(g)} }{n_\Gamma(g)} .
\end{eqnarray*}
We use the $\Psi$-construction as in \cite{bunkeolbrich94}.
Let $f(p)$ be a function of one complex variable and $
N\in\Naaa$. Then $\Psi(f(p))$
is a function of $N$ complex variables $p_1,\dots,p_N$.
It is defined if $p_i\not=p_j,\forall i\not=j$.
For example, consider the function $e^{pt}$. Then
$$\Psi(e^{pt}):=\sum_{i=1}^N\left(\prod_{j=1,j\not= i}^N
\frac{1}{p^2_j-p^2_i}e^{p_it}\right).$$
Let $R(p^2):=(A_M^2+p^2)^{-1}$ and $L<\infty$ with
 $-L<A_M^2$.
Using our distributional trace formulas we can
 (similar to \cite{bunkeolbrich94}, Prop.5.5) prove
\begin{prop}[$D$ and the trace of resolvents]
\label{traceprod} For $N\ge n/2+1$ and
$p_1,\dots,p_N\in \C$ with $Re(p_i)>L$ and $p_i\not= p_j$ for all $i\not= j$
we have in the non-super case
$$Tr\prod_{i=1}^NR(p_i^2)
=d\pi\:vol(M)\omega_{n+1}^{-1}r \Psi(\frac{P(p)}{p})+\Psi(\frac{D(p)}{2p}).$$
In the super case of (b)
we have
$$Tr\:D_M \prod_{i=1}^NR(p_i^2)
=-\imath\Psi(\frac{D^s(p)}{2}).$$
\end{prop}
{\it Proof:}
We discuss first the non-super case.
Consider the function
$$f(pt)=\frac{1}{2}\left\{\begin{array}{cc} e^{-pt}&t>0\\
                                - e^{pt}&t\le 0\end{array}\right. .$$
If $Re(p_i)>0,\forall i=1,\dots,N$, then
$\Psi(f(pt))$ vanishes exponentially if  $t\to \pm\infty$.
Moreover, it is anti-symmetric with respect to $t$ and smooth
for $t\not=0$. At $t=0$ is has continuous derivatives up to the
order $2N-1$, since the $\Psi$ kills the first $N-1$ even Taylor
coefficients \cite{bunkeolbrich94}, Lemma 5.4.

Note that
\begin{eqnarray*}
\prod_{i=1}^NR(p_i^2)&=&\int_0^\infty \Psi(e^{-pt})\frac{sin(tA_M)}{A_M} dt\\
&=&\int_{-\infty}^\infty  \Psi(f(pt)) \frac{sin(tA_M)}{A_M} dt.
\end{eqnarray*}
Integrating the distributional trace formula
for $Tr\:cos(tA_M)$ we obtain
\begin{eqnarray*}
&&Tr \frac{sin(tA_M)}{A_M}\\&=&vol(M)\omega_{n+1}^{-1}\:\pi r
P(\frac{d}{dt})(\Theta(t)-1/2)\\&&+\sum_{[g]\in C\Gamma,[g]\not=[1]}
\frac{C(g,\sigma)tr\:\chi(g)}{n_\Gamma(g)}
(\Theta(t-l(g))-\Theta(-t-l(g)))\quad \mbox{case (a)},      \\[4mm]
&&Tr \frac{sin(tA_M)}{A_M}\\&=&vol(M)
\omega_{n+1}^{-1}\:2\pi r P(\frac{d}{dt})
(\Theta(t)-1/2)\\&&+\sum_{[g]\in C
\Gamma,[g]\not=[1]} \frac{(C(g,\sigma)+C(g,w\sigma))tr\:
\chi(g)}{n_\Gamma(g)}(\Theta(t-l(g))-\Theta(-t-l(g)))\quad \mbox{case (b)},
\end{eqnarray*}
where the integration constant was fixed such that the right hand
side
is anti-symmetric.
We can insert $\Psi(f(pt))$ into the trace formula.
In order to justify this, note that one can apply the trace
formula to an anti-symmetric version $\tilde{\Psi(f(pt))}$
   smoothed
 at $t=0$ such that
$\tilde{\Psi(f(pt))}=\Psi(f(pt))$ for $|t|>1$.
Now let a sequence of such $\tilde{\Psi(f(pt))}$ tend to  $\Psi(f(pt))$
in $C^{2N-1}_{loc}(\R)$. Then
$$\int_{-\infty}^\infty  \tilde{\Psi(f(pt))}
\frac{sin(tA_M)}{A_M} dt\to
  \int_{-\infty}^\infty  \Psi(f(pt)) \frac{sin(tA_M)}{A_M} dt$$
in the sense of trace class operators.
Moreover, all individual terms in the trace formula converge
such that their sum converges, as well.
We obtain
\begin{eqnarray*} &&Tr\prod_{i=1}^NR(p_i^2)\\
&=&d\pi r\omega_{n+1}^{-1}vol(M)
\int_{-\infty}^\infty (P(\frac{d}{dt})(\Theta(t)-1/2)) \Psi(f(pt))  dt\\
&+&\sum_{[g]\in C\Gamma,[g]\not=1}
\frac{\left\{\begin{array}{cc}
C(g,\sigma)tr\:\chi(g)&\mbox{(a)}\\(C(g,\sigma)+
C(g,w\sigma))tr\:\chi(g)&\mbox{(b)}\end{array}\right\}}{n_\gamma(g)}
\Psi(\frac{f(pl(g))}{p})\end{eqnarray*}
The second term is nothing else then $$\Psi(\frac{D(p)}{2p})$$
while the first  term gives
\begin{eqnarray*}
&&2d\pi r\omega_{n+1}^{-1} vol(M) \int_{0}^\infty P(\frac{d}{dt}) \Psi(f(pt))
dt\\
&=&d\pi r\omega_{n+1}^{-1} vol(M)\int_{0}^\infty P(\frac{d}{dt})
\Psi(e^{-pt})dt\\
&=&d\pi r\omega_{n+1}^{-1} vol(M)\int_{0}^\infty \Psi(P(p)e^{-pt})dt\\
&=&d\pi r\omega_{n+1}^{-1} vol(M)\Psi(\frac{P(p)}{p})
\end{eqnarray*}

In the super case we argue similarly. We have
$$D_M \prod_{i=1}^NR(p_i^2)=\int_{-\infty}^\infty\Psi(f(pt))
D_M \frac{sin(tA_M)}{A_M} dt.$$
Applying the integrated distributional trace formula
$$Tr\:D_M  \frac{sin(tA_M)}{A_M}=-\imath\sum_{[g]\in C\Gamma,[g]\not=[1]}
\frac{C(g,\sigma)-C(g,w\sigma)}{n_\Gamma(g)}(\delta(t-l(g))-\delta(t+l(g)))$$
we obtain
\begin{eqnarray*}  Tr \:D_M\prod_{i=1}^NR(p_i^2)
&=&\sum_{[g]\in C\Gamma,[g]\not=1}
-\imath\frac{(C(g,\sigma)-C(g,w\sigma))tr\:\chi(g)}{n_\Gamma(g)}
\Psi(e^{-pl(g)})\\
&=&-\imath\Psi(\frac{D^s(p)}{2}).
\end{eqnarray*}
$\Box$\newline
\begin{kor}
The logarithmic derivative of the Selberg zeta function
$D^{(s)}(p)$ has an analytic continuation to all of $\C$.
Its singularities are first order poles at $\pm\imath \lambda$,
$\lambda\in spec\:A_M$
with residue $dim\: E_{A_M}(\{\lambda\})$ ($2dim\: E_{A_M}(\{0\})$ if
 $\lambda=0$) in the non-super case and
at $\imath\lambda$, $\lambda\in spec\: D_M$, with residue
$dim\:E_{D_M}(\{\lambda\})-dim\:E_{D_M}(\{-\lambda\})$ in the
 super case.
\end{kor}
\subsection{The fundamental properties of the zeta functions}
All residues of the logarithmic derivatives of $D(p)^{(s)}$
are integers. Hence we can
apply $exp\circ \int_{s}^\infty dp$ to $D^{(s)}(p)$ in order to
obtain the continuation of the Selberg zeta functions.
\begin{prop} \label{sing}
$Z_S(s)$ in the case (a) and S(s)
in
 the case (b) have analytic continuations to all of $\C$.
Their singularities (zeros and poles) are
at $\pm\imath \lambda$, $\lambda\in spec\:A_M$
of order $dim\: E_{A_M}(\{\lambda\})
$ ($2dim\:
 E_{A_M}(\{0\})$ if $\lambda=0$).
In the case (b) $S^s(s)$ has its singularities at
$\imath \lambda$, $\lambda\in spec\:D_M$
of order
$dim\:E_{D_M}(\{\lambda\})-dim\:E_{D_M}(\{-\lambda\})$.
\end{prop}
We derive the functional equations of the zeta functions.
\begin{prop}\label{rechnt}
The functional equations of the Selberg zeta functions are
\begin{eqnarray*}
\frac{Z_S(s)}{Z_S(-s)}&=
&exp\left(2\pi r\omega_{n+1}^{-1} vol(M)\int_{0}^s P(p) dp\right) \:\mbox{case
(a)}\\
\frac{S(s)}{S(-s)}&=&exp
\left(4\pi r \omega_{n+1}^{-1}vol(M) \int_{0}^s P(p) dp\right) \:\mbox{case
(b)}.
\end{eqnarray*}
The zeta function $S^s(s)$ satisfies the functional equation
$$S^s(s)S^s(-s)=e^{2\pi\imath\eta(D_M)}.$$
Moreover $S^s(0)=e^{\imath\pi\eta(D_M)}$.
\end{prop}
The results on $S^s$ were obtained
 in the special case of the signature operator
by Millson \cite{millson78}.
 They also follow by a specialization
of results of Moscovici/Stanton
\cite{moscovicistanton89}
 as explained in Juhl \cite{juhl93}.
{\it Proof:}
We first consider the non-super case.
{}From Proposition \ref{traceprod} we obtain by inserting $p_1=\pm p$ and
summing
\begin{eqnarray*}
D(p)+D(-p)&=&d\pi r \omega_{n+1}^{-1} vol(M) P(p)\\
D^s(p)-D^s(-p)&=&0.
\end{eqnarray*}
Integrating and exponentiating the first equation we obtain
 the desired result.
In the case (b) the functional equation for $S^s(s)$ is more
interesting since the integration constant does not dop out.
Note that $$R(p^2)=\int_0^\infty e^{-p^2t} e^{-t(A_M^s)^2} dt.$$
Hence, assuming $p_i>>0, \forall i=1,\dots,N$, we obtain
\begin{eqnarray*}
\Psi(D^s(p))&=&2\imath Tr\: D_M\prod_{i=1}^N  R(p_i^2)\\
&=&2\imath Tr \:D_M \Psi(R(p^2))\\
&=&2\imath Tr  \int_0^\infty \Psi(e^{-p^2t}) D_M e^{-t(A_M^s)^2} dt.
\end{eqnarray*}
Since $Tr \:D_Me^{-t(A_M^s)^2}$ is regular near $t=0$ (and
 vanishes there) and
$D^s(p)$ vanishes at infinity, we can  deduce
\begin{eqnarray*}
D^s(p)&=&2\imath \int_0^\infty e^{-p^2t} Tr \:D_M e^{-t(A_M^s)^2} dt\\
log\:S^s(s)&=&2\imath\int_0^\infty
\frac{\sqrt{\pi}}{2\sqrt{t}}(1-\Phi(t^{1/2}s)) Tr\:D_Me^{-t(A_M^s)^2} dt,
\end{eqnarray*}
where $\Phi(x):=\frac{2}{\sqrt{\pi}}\int_0^x e^{-x^2} dx$ is
the error integral.
Hence
\begin{eqnarray*}
log\:S^s(s)+log\:S^s(-s)&=&2\imath
\int_0^\infty \frac{\sqrt{\pi}}{\sqrt{t}} Tr\:D_M e^{-t(A_M^s)^2} dt\\
&=&2\pi\imath \eta(D_M).
\end{eqnarray*}
The functional equation follows.
$\Box$\newline
\begin{prop}
In the case (b) the zeta function $Z_S(s):=Z_S(s,\sigma)$
has a meromorphic continuation
to the whole complex plane. It has singularities at $0\not=s=
\imath \lambda$
of order
$$\frac{1}{2}(dim\:E_{A_M}(\{\lambda\})+dim\:E_{D_M}
(\{\lambda\})-dim\:E_{D_M}(\{-\lambda\}))$$ and
of order $dim\:ker\:A_M$ at $s=0$. $Z_S(s)$
satisfies the functional equation
$$\frac{Z_S(s,\sigma)}{Z_S(-s,w\sigma)}=
e^{\pi\imath\eta(D_M)}exp\left(2\pi r \omega_{n+1}^{-1} vol(M)
\int_{0}^s P(p) dp\right).$$
\end{prop}
{\it Proof:}
For $Re(s)>\rho$ we have $Z_S(s)=\sqrt{S(s)S^s(s)}$.
The continuation of $Z_S(s)$ is obtained from
that of $S^{(s)}(s)$ if we can take the square root. It is sufficient
to show that all singularities of $S(s)S^s(s)$ have even order.
We choose the representative $\gamma$
of the lift of $\sigma+w\sigma$ such that
$V(\gamma)=V(\gamma^s)$ if one forgets the $\Z_2$-grading.
In particular, $\lambda\in spec\:|D_M|$ iff $\lambda\in spec\:A_M$.
For $\lambda > 0$ we have
$dim\:E_{D_M}(
\{\lambda\})-
dim\:E_{D_M}(\{-\lambda\})=
dim\:E_{A_M}(\{\lambda\})$
(mod $2$).
Hence the order of the singularity at $s\not=0$ of $S(s)S^s(s)$ is even
by Proposition \ref{sing}.
At $s=0$ the function $S^s(s)$ is regular and $S(s)$ has a singularity
 of even order.
Thus we can take the square root above and obtain a
meromorphic continuation of
$Z_S(s)$ to all of $\C$.
The order singularity of $Z_S(s)$ at $0\not=s=\imath
 \lambda$
is $\frac{1}{2}(dim\:E_{A_M}(\{\lambda\})+
dim\:E_{D_M}(\{\lambda\})-dim\:E_{D_M}(\{-\lambda\}))$
and $dim\:ker\:A_M$ at $s=0$.
We also obtain
\begin{eqnarray*}
\frac{Z_S(s,\sigma)}{Z_S(-s,w\sigma)}&=
&\sqrt{\frac{S(s,\sigma)S^s(s,\sigma)}{S(-s,w\sigma)S^s(-s,w\sigma)}}\\
&=&\sqrt{\frac{S(s,\sigma)}{S(-s,\sigma)}S^s(s,\sigma)S^s(-s,\sigma)}\\
&=&e^{\pi\imath\eta(D_M)}exp\left(2\pi r \omega_{n+1}^{-1} vol(M) \int_{0}^s
P(p) dp\right)
\end{eqnarray*}
$\Box$\newline

Now we relate the zeta functions to zeta-regularized determinants.
\begin{prop}\label{detdet}
We have
\begin{eqnarray*}
Z_S(s)&=&det(s^2+A_M^2) e^{2\pi r\omega_{n+1}^{-1} vol(M) \int_{0}^s P(p) dp}\\
S(s)&=& det(s^2+A_M^2) e^{4\pi r\omega_{n+1}^{-1} vol(M) \int_{0}^s P(p) dp}.
\end{eqnarray*}
\end{prop}
{\it Proof:}
For $|Re(p_i)|$ large we consider the function $L_M$ of $p_i$, $i=1,\dots,N$ :
$$L_M(p_1,\dots,p_N):=Tr\:\Psi(\frac{1}{p^2+A^2_M}).$$
We can write
\begin{eqnarray*}L_M(p_1,\dots,p_N) &=&Tr\Psi(\int_0^\infty
e^{-t(p^2+A^2_M)}dt)\\
&=&Tr\Psi(\int_0^\infty -\frac{d}{2pdp}e^{-t(p^2+A^2_M)}\frac{dt}{t}) \\
&=&\Psi(-\frac{d}{2pdp} P.F.\int_0^\infty Tr e^{-t(p^2+A^2_M)}\frac{dt}{t}),
\end{eqnarray*}
where $P.F. \int_0^\infty\dots$ stands for taking
 the finite part of $\int_\epsilon^\infty\dots $.

Applying Soule/Abramovich/Burnol/Kramer,
\cite{souleabramovichburnolkramer92}, Theorem 5.1.1, we obtain
\begin{equation}\label{ext333}L_M(p_1,\dots,p_N)
=\Psi(\frac{d}{2pdp}ln(det(p^2+A^2_M))).\end{equation}
The determinant is the zeta regularized (super)
determinant of $p^2+A^2_M$ defined by
$$ln(det(p^2+A^2_M)):=- P.F.\int_0^\infty Tr e^{-t(p^2+A^2_M)} \frac{dt}{t}$$
(see the explanation on p.39 in \cite{bunkeolbrich94}).
Equation (\ref{ext333}) extends analytically  to all of $\C$.
We obtain
$$\Psi(\frac{D(p)+2d\pi r\omega_{n+1}^{-1}
vol(M)P(p)}{2p})=\Psi(\frac{d}{2pdp} ln(det(p^2+A^2_M))).$$
It follows
$$D(p)+2d\pi r\omega_{n+1}^{-1} vol(M)P(p) =
R^\prime(p)+\frac{d}{dp} ln(det(A^2_M+p^2)),$$
where $R^\prime(p)$ is a certain odd polynomial.
Integrating once and exponentiating we obtain
that the Selberg zeta function has the representation
$$
Z(s)=e^{R(s)}det(s^2+A^2_M)e^{d\pi r
\omega_{n+1}^{-1} vol(M) \int_{-s}^s P(p) dp},$$
where $Z(s)$ is $Z_S(s)$ ($d=1$) or $S(s)$ ($d=2$), respectively.
\begin{lem}
$R(s)=0$
\end{lem}
{\it Proof:}
Note that $log \:D(s)\to 0$ exponentially as $s\to 0$.
Let
$$Tr e^{-tA_M^2}\stackrel{t\to 0}{\sim}
\sum_{k=0}^\infty c_{(2k-n)/2} t^{(2k-n)/2}$$
define the real numbers $c_{(2k-n)/2}$.
Then by \cite{voros87}, Eq. 5.1,
$$-log\:det(s^2+A^2_M)\stackrel{s\to \infty}{\sim}
\sum_{k=0}^\infty  c_{(2k-n)/2} \Gamma((2k-n)/2) s^{n-2k}.$$
This asymptotic expansion containing only odd powers of $s$
must cancel with the Plancherel term, i.e.
$$d\pi vol(M) \int_{-s}^s P(p) dp =  \sum_{k=0}^{(n-1)/2}  c_{(2k-n)/2}
\Gamma((2k-n)/2) s^{n-2k},$$
$c_{(2k-n)/2}=0$ for $2k>n$.
It also follows $R(p)=0$.
$\Box$\newline
\section{The Ruelle zeta function}\label{ruelleee}

\subsection{Definition and relation with the Selberg zeta function}

We discuss the Ruelle zeta function of compact odd dimensional
hyperbolic manifolds $M$.
Ruelle \cite{ruelle76} introduced a zeta function for
Anosov flows coding the length spectrum of closed orbits.
In our situation the corresponding flow is the geodesic flow
on the unit sphere bundle of $M$ (twisted with a flat vector bundle).
The closed orbits correspond
to closed geodesics and thus to the conjugacy classes of $\Gamma$.
\begin{ddd}
The Ruelle zeta function of $M$ associated to a r-dimensional unitary
representation $\chi$ of $\Gamma$
is defined by the infinite product
$$Z_{R,\chi}(s):=\prod_{[g]\in C\Gamma,[g]\not=1,primitive}
det(1-\chi(g)e^{-sl(g)})^{-1}$$
converging for $Re(s)>2\rho$.
\end{ddd}
We consider the representation $\chi$ as fixed throughout this and the
following subsection and shall omit the corresponding index.

Let $\sigma^p$ be the p-th exterior powers of the
complexified standard representation of $SO(n-1)$.
For $p\not=(n-1)/2$ the representation $\sigma^p$ is
irreducible, while $\sigma^\frac{n-1}{2}$
splits into two
 non-equivalent irreducible subrepresentations $\sigma^\frac{n-1}{2}=
\sigma^+\oplus\sigma^-$, the spaces of selfdual and anti-selfdual forms.
We have for the non trivial element $w$ of the Weyl group $w\sigma^+
=\sigma^-$. Therefore it is natural to set
\[ Z_S(s,\sigma^\frac{n-1}{2}):= S(s,\sigma^+)\ .\]
The following proposition is due to Fried (\cite{fried862},\cite{fried86}). One
may also compare this with the computations for the even dimensional case in
\cite{bunkeolbrich94}, 7.1, 8.1.
\begin{prop}\label{prewq}
The Ruelle zeta function has the representation
\begin{eqnarray*}
Z_R(s)&=&\prod_{p=0}^{n-1} Z_S(s+\frac{n-1}{2}-p,\sigma^p)^{(-1)^p}\\
&=&S(s,\sigma^+)^{(-1)^\frac{n-1}{2}}\prod_{p=0}^{\frac{n-3}{2}}
\left(Z_S(s+\frac{n-1}{2}-p,\sigma^p)
Z_S(s-(\frac{n-1}{2}-p),
\sigma^p)
\right)^{(-1)^p}\ .
\end{eqnarray*}
\end{prop}
For the second equation one uses Poincar\'e duality.

The representation of the Ruelle zeta function provides
its meromorphic continuation to the whole complex plane.
That the Ruelle zeta function has a meromorphic
continuation was observed in Ruelle's original work by dynamical methods
(\cite{ruelle76}).

\subsection{The functional equation for the Ruelle zeta function}

We employ the functional equations for the Selberg zeta functions
in order to provide a functional equation for
 the
 Ruelle zeta function. We obtain
\begin{eqnarray}
\frac{Z_R(s)}{Z_R(-s)}&=&\left(
\frac{S(s,\sigma^+)}{
S(-s,\sigma^+)}\right)^{(-1)^\frac{n-1}{2}}
\prod_{p=0}^\frac{n-3}{2}
\left(\frac{Z_S(s+\frac{n-1}{2}-p,\sigma^p)Z_S(s-(\frac{n-1}{2}-p),\sigma^p)}
{Z_S(-s+\frac{n-1}{2}-p,\sigma^p)Z_S(-s-(\frac{n-1}{2}-p),\sigma^p)}
\right)^{(-1)^p}\nonumber\\
&=&exp
\left(
\frac{2\pi r vol(M)}{\omega_{n+1}}\ \big( 2(-1)^\frac{n-1}{2}\int_0^s
P(q,\sigma^+)dq+\right.\label{fd}\\
 &&\left.\qquad +\sum_{p=0}^{\frac{n-3}{2}} (-1)^p (
\int_{0}^{s+\frac{n-1}{2}-p} P(q,\sigma^p)dq+\int_{0}^{s-(\frac{n-1}{2}-p)}
P(q,\sigma^p)dq)\big)
\right)
\ .\nonumber
\end{eqnarray}
Let $h(s)$ denote the derivative of the exponent
(up to the prefactor), i.e. the polynomial
\begin{eqnarray*}
h(s)=2(-1)^\frac{n-1}{2} P(s,\sigma^+)+\sum_{p=0}^{\frac{n-3}{2}} (-1)^p (
P(s+\frac{n-1}{2}-p,\sigma^p)+P(s-(\frac{n-1}{2}-p),\sigma^p))\ .
\end{eqnarray*}
\begin{lem}
The polynomial $h(s)$ is a constant
$$h(s)\equiv n+1\ .$$
\end{lem}
{\it Proof:}
The claim is a consequence of identities in the Weyl polynomials.
In order to establish them we employ a geometric argument.
It is enough to evaluate $h(s)$ at a large number of points.

We consider the Grassmannian $B=SO(n+1)/SO(n-1)\times
SO(2)$ of oriented 2-planes in $\R^{n+1}$ as a homogeneous K\"ahler manifold.
Each pair
$(\sigma,k)\in SO(n-1)^\wedge  \times \Z=(SO(n-1)\times SO(2))^\wedge $ defines
a homogeneous
holomorphic vector bundle $W(k,\sigma)$ on $B$.

The Borel-Weil-Bott theorem asserts that the representation
of $SO(n+1)$ with the highest weight $\Lambda=
\mu_\sigma+k\alpha\in\imath \taaa^\ast\oplus   \aaa^\ast$  ($k>>0$)
can be realized as the space of holomorphic sections of $W(k,\sigma)$ and
all higher cohomology groups of $W(k,\sigma)$ vanish.
A consequence of the theorem of Borel-Weil-Bott is
$$P(k+\frac{n-1}{2},\sigma)=\chi_a(B,W(k,\sigma)), \quad k\in\Z,$$
 where $\chi_a$ is the analytic genus of the bundle $W(k,\sigma)$,
i.e. the
Euler characteristic of the complex given by the Dolbeault resolution
of $W(k,\sigma)$.
Thus
$$P(k+\frac{n-1}{2},\sigma)=index(\bar{\partial}+\bar{\partial}^\ast),$$
i.e.
$P(k+\frac{n-1}{2},\sigma)$ is the index of the Dirac type operator
$\bar{\partial}+\bar{\partial}^\ast$ on the $\Z_2$-graded vector bundle
$\Lambda^{0,\ast}T^\ast B\otimes W(k,\sigma)$
(the grading given by even and odd form degree).

Note that $\Lambda^{p,0}T^*B\cong W(-p,\sigma^p)$. We obtain for
 $s\in \Z$ and $p\leq \frac{n-3}{2}$
\[ P(s+\frac{n-1}{2}-p,\sigma^p)=\chi_a(\Lambda^{p,0}T^*B\otimes
W(s,1))\]
and
\[ P(s-(\frac{n-1}{2}-p),\sigma^p)=\chi_a(\Lambda^{n-1-p,0}T^*B\otimes W(s,1))\
.\]
In addition,
\[2 P(s,\sigma^+)=P(s,\sigma^+)+P(s,\sigma^-)=
\chi_a(\Lambda^{\frac{n-1}{2},0}T^*
B\otimes W(s,1))\ .\]
It follows that $h(s),s\in\Z$ is the
index of $(\bar{\partial}+\bar{\partial}^\ast)$
on the bundle
$$\Lambda^{\ast,\ast}T^\ast B\otimes W(s ,1)=\Lambda^\ast T^\ast
 B\otimes W(s ,1)$$
graded by the total form degree.
Now the differential operator $(\bar{\partial}+\bar{\partial}^\ast)$
can be deformed to $D=\nabla+\nabla^\ast$, where $\nabla$ is induced by the
Levi-Civita connection of $B$ and the homogeneous connection on $W(s+\rho,1)$.
Here $\nabla$ acts by alternating differentiation.
The deformation is given by
$$D_t=(\bar{\partial}+\bar{\partial}^\ast)+t(\partial+\partial^\ast)$$
and stays inside the elliptic operators.
$\partial$ is again defined using the connection.
The index of $D_t$ is independent of $t$ and for $t=1$ it is independent
of the twisting line bundle and equal to $index (D_1)=\chi(B)$.
The Euler characteristic of $B$ can be computed using the
orders of Weyl groups (see for example Bott \cite{bott65})
\[ \chi(B)=\frac{|W(SO(n+1))|}{|W(SO(n-1)\times SO(2))|}=
\frac{|W(SO(n+1))|}{|W(SO(n-1))|}=\frac{2^\frac{n-1}{2} (\frac{n+1}{2})!}
{2^\frac{n-3}{2} (\frac{n-1}{2})!}=n+1\ .\Box\]

Now we come back to discuss the functional equation of the
 Ruelle zeta function. $h(s)=n+1$ is the derivative of a polynomial
$(n+1)s+C$ for some constant $C$. But the left hand side of (\ref{fd})
evaluated at $s=0$ is $1$. We conclude that $C=0$ and obtain
\begin{theorem}
The Ruelle zeta function satisfies the functional equation
$$\frac{Z_R(s)}{Z_R(-s)}=  e^{\frac{2\pi r(n+1)vol(M)}{\omega_{n+1}} s}.$$
\end{theorem}

\subsection{Analytic torsion and the Ruelle zeta function}

We derive the result of Fried \cite{fried862}
that was generalized by Moscovici/Stanton
\cite{moscovicistanton91} to higher rank situations.
Assume that the twist is acyclic, i.e. $H^\ast(\Gamma,\chi)=0$.
Recall that the analytic torsion   of $M$ with respect to the twist $\chi$
is then defined by
$$\tau=\tau_M(\chi):=\sqrt{\prod_{l=1}^n (det\Delta_{\chi,l})^{(-1)^ll}},$$
where $\Delta_{\chi,l}$ is the Laplacian on $l$-forms on $M$ twisted with
the flat bundle associated to $\chi$. In the following we will omit $\chi$ in
our
notation.
\begin{theorem}
$Z_R(0)=\tau^{-2}$.
\end{theorem}
{\it Proof:}
Using Poincare  duality, Proposition \ref{prewq} and Proposition \ref{detdet}
we obtain
\begin{eqnarray*}
Z_R(0)&=&\prod_{p=0}^{n-1} Z_S(\frac{n-1}{2}-p,\sigma^p)^{(-1)^p}\\
&=&\prod_{p=0}^{n-1} det(A_M^2(\sigma^p)+(\frac{n-1}{2}-p)^2)^{(-1)^p}\\
&=&det(A_M(\sigma^++\sigma^-))^{(-1)^{(n-1)/2}}\prod_{p=0}^{(n-3)/2}
det(A_M^2(\sigma^p)+(\frac{n-1}{2}-p)^2)^{2(-1)^p}.
\end{eqnarray*}
Now $\sigma^p=r(\sum_{l=0}^p (-1)^{l-p} \lambda^l)$ for $p\le (n-1)/2$ and
$\Delta_l=A_M^2(\sigma^l)+(\frac{n-1}{2}-l)^2 $.
Thus
\begin{eqnarray*}Z_R(0)&=&\left(
\prod_{l=0}^{(n-1)/2}det(\Delta_l)^{(-1)^{l}}\right)
\prod_{p=0}^{(n-3)/2}
\prod_{l=0}^{p}
det(\Delta_l)^{2(-1)^{l}}\\&=&
\prod_{l=0}^{(n-1)/2} det(\Delta_l)^{2(\frac{n}{2}-l)(-1)^l}\\
&=&\tau^{-2}.
\end{eqnarray*}
\bibliographystyle{plain}


\end{document}